\documentclass[a4paper,10pt]{article}
\usepackage[cp1251]{inputenc}
\usepackage[english]{babel}

\usepackage{amssymb}\usepackage{amsmath}
\usepackage{epsfig}
\usepackage{psfig}
\def\bl{\rule[-1mm]{2.4mm}{2.4mm}}

\def\be{\begin{equation}}
\def\ee{\end{equation}}
\newtheorem{theorem}{\bf Theorem}
\newtheorem{lemma}{\bf Lemma}

\begin{document}

\title {Elementary construction of some Jenkins-Strebel differentials} 
\author{\copyright 2010 ~~~~A.B.Bogatyrev
\thanks{Partially supported by RFBR grants 10-01-00407 and RAS Program 
"Modern problems of theoretical mathematics"}} 
\date{} 
\maketitle

The behavior of trajectories of a regular quadratic differential on a Riemann surface is chaotic but in one case. When the critical graph of the foliation is compact, its complement is a finite set of cylinders swept out by homotopic closed leaves. The existence of such differentials was proved by Jenkins and Strebel yet in 1950-1960 \cite{J,S1,S2}. Later a more simple proof was given by Wolf \cite{WM}. Mentioned results are purely existence theorems and there are very few explicit constructions of those differentials. See e.g. \cite{A,ALS} for one- parametric families of them. Of course, one should not think that Jenkins-Strebel (or JS) differentials are somewhat exceptional. After all, they are dense in the space of regular quadratic differentials. The purpose of this note is to give an explicit multi-parametric construction for JS differentials on real algebraic curves. Roughly speaking, the square of any real holomorphic abelian differential subjected to certain linear restrictions will be a JS quadratic differential.

Let a compact genus $g$ Riemann surface $X$ admits an anticonformal involution  $\bar{J}$ which we also call reflection. The components of the set of fixed points of the involution are smooth closed curves known as real ovals \cite{N}. The reflection $\bar{J}$ acts naturally on the space of 1-homologies of the surface and splits it into the subspaces corresponding to eigenvalues $\pm1$ of operator $\bar{J}$:

\be
\label{H1split}
\mathbb{R}^{2g}\cong
H_1(X, \mathbb{R})=
H_1^+(X, \mathbb{R})\oplus
H_1^-(X, \mathbb{R}),
\qquad
H_1^\pm(X, \mathbb{R}):=(I\pm\bar{J})H_1(X, \mathbb{R}).
\ee

The elements $C=\bar{J}C$ of the subspace  $H_1^+(X)$ we call  \emph{even cycles}.
Respectively,  the elements $C=-\bar{J}C$ of the subspace $H_1^-(X)$ we call \emph{odd cycles}.
The subspaces of even and odd cycles contain the full rank lattices of integer cycles
$H_1^\pm(X,\mathbb{Z}):=H_1^\pm(M,\mathbb{R})\cap H_1(X,\mathbb{Z})$.
For example, the oriented real ovals of the surface (if any) are its even integer cycles.

The intersection index  of cycles is a non-degenerate skew-symmetric form on the space
$H_1(M, \mathbb{R})$. The reflection $\bar{J}$ changes the orientation in any intersection point of integer cycles,
therefore
\be
\bar{J}C_1\circ\bar{J}C_2=-C_1\circ C_2,
\qquad C_1,C_2\in H_1(M, \mathbb{R}).
\ee
Now we easily see that the subspaces of even and odd cycles are lagrangian (i.e.  the intersection form vanishes  in the subspaces), their dimensions are equal and therefore equal to $g$.

The space $\Omega^1(X)\cong\mathbb{C}^g$ of holomorphic differentials on $X$ contains the subspace of the  
so called real differentials $\Omega^1_\mathbb{R}(X)\cong\mathbb{R}^g$,
which become complex conjugate after the action of reflection:
$\bar{J}^*\eta:=\overline{\eta}$, $\eta\in\Omega^1_\mathbb{R}(X) $.
The integration of real differentials along even (resp. odd) cycles
give us real  (resp. pure imaginary) numbers:
$$
\int_C\xi=
\int_{\pm\bar{J}C}\xi=
\pm\int_C\bar{J}^*\xi=
\pm\int_C\overline{\xi}=
\pm\overline{\int_C\xi}.
$$

\begin{lemma}
The space $(H_1^-(X, \mathbb{R}))^*\cong\mathbb{R}^g$ of real linear functionals over odd cycles is canonically isomorphic to each of the following two spaces:\\
(i) $H_1^+(X, \mathbb{R})$;  \hskip 2cm 
(ii) $\Omega^1_\mathbb{R}(X)$. 
\end{lemma}
(i). In this case the functional is given by the intersection form. The non-degeneracy of this form implies
that an even cycle annihilating all odd cycles should be zero.\\
(ii). In this case the functional is given by the formula
$$
\langle \eta|C^-\rangle:=i\int_{C^-}\eta,
\qquad C^-\in H_1^-(X). 
$$
If real differential $\eta$ annihilates all odd cycles, then all of its periods are real and therefore $\eta=0$. ~~\bl

{\bf Remarks}:
1) Usually, for the normalization of abelian differentials they take half of a canonical basis in the homologies: 
$A$- cycles or $B$-cycles. From statement (ii) of the above lemma it follows that one can use either even or odd cycles for this purpose once the surface admits a reflection. To generalize this observation let us show that any $g$-dimensional lagrangian subspace of the homology space may be used for the normalization. In other words, \emph{there exist a unique holomorphic differential on the surface $X$ with given periods in the basis of the lagrangian subspace of homologies}.\\
We choose the basis  $C_1,C_2,\dots,C_{2g}$ in the space of real homologies of the curve $X$ so that
the first  $g$ elements are in the given lagrangian subspace. We do not assume that this basis is either canonical or even integer. Riemann bilinear identity holds: 
\be
\label{volume}
0\le||\eta||^2=i\int_X
\eta \wedge \overline{\eta}= -i\sum\limits_{s,j=1}^{2g}F_{sj}\int_{C_s}\eta\overline{\int_{C_j}\eta},
\ee
where the matrix $F_{sj}$ is the inverse of the intersection matrix $C_s\circ C_j$. 
If $\int_{C_j}\eta=0$ for $j=1,\dots,g$, then the sum in the right-hand side contains the terms with $s,j>g$ only.  But in the latter case $F_{sj}=0$. Indeed, the intersection matrix has block $2\times2$ structure with zero 
$g\times g$ matrix in (1,1) position. The inverse matrix has zero block of the same size in (2,2) position. We see 
that only zero holomorphic differential has zero periods along all cycles of our lagrangian subspace.

2)From the above lemma it follows that there is a 1-1 correspondence between  even cycles  $C^+$ and real differentials $\eta$ on the surface given by the rule:
$$
i\int_C\eta=C^+\circ C,
\qquad\forall C\in H_1^-(X,\mathbb{R}).
$$
This correspondence is not Poincare duality as the latter equality is not true for even cycles $C$.

3) Let the surface $X$ has $k$ real ovals. The span of real ovals in the homology space has real dimension $k$
if the surface is not \emph{separating} (i.e. the surface with  real ovals removed is connected). Otherwise the dimension is 1 less. We'll show that the differentials corresponding to the points of this linear span are JS when squared. Therefore we give a $k$ or $(k-1)$ -parametric family of JS differentials on a given surface $X$.

\begin{figure}[h!]
\begin{picture}(150,55)
\put(50,2){\psfig{figure=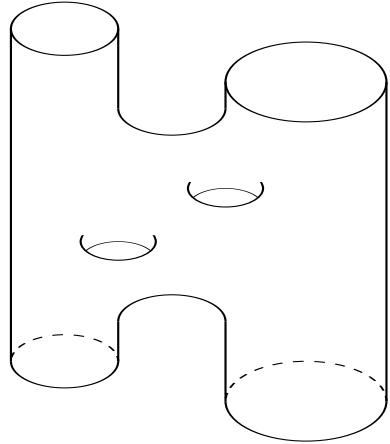}}
\end{picture}
\caption{\small Function $H(x)$ as a height function on the surface with real ovals removed.}
\label{H}
\end{figure}

\begin{theorem}
Let the integral of real holomorphic differential $\eta$ along an odd cycle $C^-$ vanishes if the intersection index of  $C^-$ with any real oval is zero. Then the foliation  $\eta^2>0$ is Jenkens-Strebel.
\end{theorem}

Proof.
Let us remove real ovals from the surface. In the remaining part(s) of the surface one can correctly define the function
$$
H(x):=Im~\int_*^x\eta,
\qquad
x\in X\setminus \{ \mbox{\rm real ovals}\}.
$$
Indeed, if the closed path $C$ does not intersect real ovals, then
$$
2Im~\int_C\eta=Im~\int_{C-\bar{J}C}\eta=0,
$$
because $(C-\bar{J}C)\circ C^+=2C\circ C^+=0$, where $C^+$ is any real oval.

This globally defined function $H(x)$ is constant on the boundary components of the cut surface and its level lines are the leaves of the foliation $\eta^2>0$ -- see Fig. \ref{H}. ~~~\bl

{\bf Example}. Let us consider a hyperelliptic curve with real branchpoints. The span of real ovals is the entire space of even cycles. Therefore, the square of any real holomorphic differential is JS.
For general real hyperelliptic curves the bases in the lattices of even and odd integer cycles are described in \cite{B1}. This makes the condition of Theorem 1 being  explicit.

\vspace{5mm}
\parbox{9cm}
{\it
119991 Russia, Moscow GSP-1, ul. Gubkina 8,\\
Institute for Numerical Mathematics,\\
Russian Academy of Sciences\\[3mm]
{\tt gourmet@inm.ras.ru, ab.bogatyrev@gmail.com}}

\end{document}